\documentclass[11pt]{article}
\usepackage{amsfonts}
\usepackage{amssymb}
\usepackage{mathrsfs}
\usepackage{float}
\usepackage{amsmath}
\usepackage{amsfonts,amssymb}
\usepackage{graphicx,color}
\catcode`\@=11 \@addtoreset{equation}{section}

\catcode`\@=12

\newtheorem{Theorem}{Theorem}[section]
\newtheorem{Definition}{Definition}[section]
\newtheorem{Lemma}{Lemma}[section]
\newtheorem{Example}{Example}[section]
\newtheorem{Corollary}{Corollary}[section]
\newtheorem{Remark}{Remark}[section]

\def\2{{I \hskip -1.0mm I}}
\def\3{{I \hskip -1.0mm I\hskip -1.0mm I}}
\def\4{{I \hskip -0.9mm V}}
\def\6{{V \hskip -1.35mm I}}

\begin{document}

\title{The hyperbolic geometric flow on Riemann surfaces}
\author{De-Xing Kong\thanks{Center of Mathematical Sciences, Zhejiang University,
Hangzhou 310027, China;}, $\quad$ Kefeng Liu\thanks{Department of Mathematics,
UCLA, CA 90095, USA;} $\quad$ and $\quad$ De-Liang Xu\thanks{Department of
Mathematics, Shanghai Jiao Tong University, Shanghai 200240, China.}}
\date{}
\maketitle

\begin{abstract}
In this paper the authors study the hyperbolic geometric flow on
Riemann surfaces. This new nonlinear geometric evolution equation
was recently introduced by the first two authors motivated by
Einstein equation and Hamilton's Ricci flow. We prove that, for any
given initial metric on ${\mathbb{R}}^{2}$ in certain class of
metrics, one can always choose suitable initial velocity symmetric
tensor such that the solution exists for all time, and the scalar
curvature corresponding to the solution metric $g_{ij}$ keeps
uniformly bounded for all time; moreover, if the initial velocity
tensor is ``large" enough, then the solution metric $g_{ij}$
converges to the flat metric at an algebraic rate. If the initial
velocity tensor does not satisfy the condition, then the solution
blows up at a finite time, and the scalar curvature $R(t,x)$ goes to
positive infinity as $(t,x)$ tends to the blowup points, and a flow
with surgery has to be considered. The authors attempt to show that,
comparing to Ricci flow, the hyperbolic geometric flow has the
following advantage: the surgery technique may be replaced by
choosing suitable initial velocity tensor. Some geometric properties
of hyperbolic geometric flow on general open and closed Riemann
surfaces are also discussed. \vskip 6mm


\vskip 3mm


\end{abstract}

\newpage\baselineskip=7mm

\section{Introduction}

Let $\mathscr{M}$  be an $n$-dimensional complete Riemannian
manifold with Riemannian metric $g_{ij}$. The following general
evolution equation for the metric $g_{ij}$
$$\frac{\partial^{2}g_{ij}}{\partial
t^{2}}+2R_{ij}+\mathscr{F}_{ij}\left(g,\frac{\partial g}{\partial
t}\right)=0\eqno{(HGF)}$$ has been recently introduced by Kong and
Liu \cite{kl} and named as {\it general version of hyperbolic
geometric flow}, where $\mathscr{F}_{ij}$ are some given smooth
functions of the Riemannian metric $g$ and its first order
derivative with respect to $t$. The most important three special
cases are the so-called {\it standard hyperbolic geometric flow} or
simply called {\it hyperbolic geometric flow}, the {\it Einstein's
hyperbolic geometric flow} (see (1.1) and (1.12) below,
respectively) and the {\it dissipative hyperbolic geometric flow}
(see \cite{k2} or \cite{dkl2}). The present paper concerns the first
two cases on Riemann surfaces.

In this paper we mainly study the evolution of a Riemannian metric
$g_{ij}$ on a Riemann surface $\mathscr{M}$ by its Ricci curvature
tensor $R_{ij}$ under the \textit{hyperbolic geometric flow}
equation
\begin{equation}
\label{1.1}\dfrac{\partial^{2}g_{ij}}{\partial t^{2}}=-2R_{ij}.
\end{equation}
As the first step of our research on this topic, we are interested
in the following initial metric on a surface of topological type
${\mathbb{R}}^{2}$
\begin{equation}
\label{1.2}t=0:~~~ds^{2}=u_{0}(x)(dx^{2}+dy^{2}) ,
\end{equation}
where $u_{0}(x)$ is a smooth function with bounded $C^{2}$ norm and satisfies
\begin{equation}
\label{1.2a}0<m\leqslant u_{0}(x)\leqslant M<\infty,
\end{equation}
in which $m, M$ are two positive constants. we shall prove the following result.

\begin{Theorem}
Given the initial metric (\ref{1.2}) with (\ref{1.2a}), for any
smooth function $u_{1}(x)$ satisfying

(a)~~$u_{1}(x)$ has bounded $C^{1}$ norm;

(b)~~it holds that
\begin{equation}
\label{1.3}u_{1}(x)\geqslant\frac{|u_{0}^{\prime}(x)|}{\sqrt{u_{0}(x)}%
},~~\forall\; x\in\mathbb{R},
\end{equation}
the Cauchy problem
\begin{equation}%
\begin{cases}
& \dfrac{\partial^{2}g_{ij}}{\partial t^{2}}=-2R_{ij}~~(i,j=1,2),\vspace
{2mm}\\
& \label{1.4} t=0:~~g_{ij}=u_{0}(x)\delta_{ij},~~\dfrac{\partial g_{ij}%
}{\partial t}=u_{1}(x)\delta_{ij}~~(i,j=1,2)
\end{cases}
\end{equation}
has a unique smooth solution for all time $t\in\mathbb{R}$, and the
solution metric $g_{ij}$ possesses the following form
\begin{equation}
\label{1.5}g_{ij}=u(t,x)\delta_{ij}
\qquad(i,j=1,2).\qquad\qquad\square
\end{equation}
\end{Theorem}

Theorem 1.1 will be proved in Section 2. This theorem gives a global
existence result on smooth solutions of hyperbolic geometric flow.
Based on Theorem 1.1 we can further prove the following theorem in
Section 3.
\begin{Theorem}  {\rm (I)}
Under the assumptions mentioned in Theorem 1.1, the Cauchy problem
(1.5) has a unique smooth solution with the form (1.6) for all time,
moreover the scalar curvature $R(t,x)$ corresponding to the solution
metric $g_{ij}$ remains uniformly bounded, i.e.,
\begin{equation}
\label{1.3}|R(t,x)|\le k,\quad\forall\;
(t,x)\in\mathbb{R}^+\times\mathbb{R},
\end{equation}
where $k$ is a positive constant depending on $M$, the $C^{2}$ norm
of $u_{0}$ and $C^{1}$ norm of $u_{1}$, but independent of $t$ and
$x$.

{\rm (I\!I)} Under the assumptions mentioned in Theorem 1.1, suppose
furthermore that there exists a positive constant $\varepsilon$ such
that
$$u_{1}(x)\geqslant\frac{|u_{0}^{\prime}(x)|}{\sqrt{u_{0}(x)}%
}+\varepsilon,~~\forall\; x\in\mathbb{R},\eqno{(1.4a)}$$ then the
Cauchy problem (1.5) has a unique smooth solution with the form
(1.6) for all time, moreover the solution metric $g_{ij}$ converges
to one of flat curvature at an algebraic rate $\dfrac
{1}{(1+t)^{\gamma}}$, i.e., $$|R(t,x)|\le
\dfrac{\tilde{k}}{(1+t)^{\gamma}},\quad\forall\;
(t,x)\in\mathbb{R}^+\times\mathbb{R},\eqno{(1.7a)}$$ where
$\gamma\in (1/2,2]$ and $\tilde{k}$ are two positive constants
depending on $\varepsilon, M$, the $C^{2}$ norm of $u_{0}$ and
$C^{1}$ norm of $u_{1}$. $\qquad\qquad\square$
\end{Theorem}

The condition (1.4) is a sufficient condition guaranteeing the
global existence of smooth solution to the Cauchy problem (1.5). On
the other hand, in some sense it is also a necessary condition,
because we have the following theorem.

\begin{Theorem} Suppose that $u_0(x)\not\equiv 0$, without
loss of generality, we may assume that there exists a point
$x_{0}\in\mathbb{R}$ such that
\begin{equation}
u_{0}^{\prime}(x_{0})<0.
\end{equation}
For the following initial velocity
\begin{equation}
u_{1}(x)\equiv\frac{u_{0}^{\prime}(x)}{\sqrt{u_{0}(x)}}%
,~~~\forall\; x\in\mathbb{R},
\end{equation}
the Cauchy problem (\ref{1.4}) has a unique smooth solution only in
$[0,\tilde{T}_{\max })\times\mathbb{R}$, moreover there exists some
point $(\tilde{T}_{\max},x_{\ast})$ such that the scalar curvature
$R(t,x)$ satisfies
\begin{equation}
R(t,x)\rightarrow\infty\qquad\text{as}~~(t,x)\nearrow (\tilde
{T}_{\max},x_{\ast}),
\end{equation}
where
\begin{equation}
\tilde{T}_{\max}=-2\left(\inf_{x\in\mathbb{R}}\left\{u_0^{\prime}(x)u_0^{-\frac32}(x)\right\}
\right)^{-1}. \quad\quad\quad\square
\end{equation}
\end{Theorem}

Theorem 1.3 will be proved in Section 4.

Motivated by Theorems 1.1-1.2, we conjecture that any complete
metric on a simply connected non-compact surface converges to a flat
metric by choosing a suitable initial velocity tensor
$\dfrac{\partial g_{ij}}{\partial t}(0,x).$ This should be true for
higher dimensional manifolds of topological type ${\mathbb R}^n$
with suitable curvature assumption. We are now working on the
general case in $\mathbb{R}^{2}$ other than (\ref{1.2})\footnote{In
fact, for the case of initial data $u_0=u_0(ax+by),u_1=u_1(ax+by)$,
by the same way we can prove some results similar to Theorems
1.1-1.3, where $a$ and $b$ are two constants satisfying $a^2+b^2\neq
1$.}.

On the other hand, Theorem 1.3 shows that if
we do not choose suitable initial velocity tensor $\dfrac{\partial g_{ij}%
}{\partial t}(0,x),$ the solution to the Cauchy problem (\ref{1.4})
blows up in finite time, and the curvature tends to infinity when
the points approach the blowup points. In this case, a flow with
surgery has to be considered. In other words, this paper attempts to
show that, by choosing a suitable initial velocity tensor, any
complete metric on a surface of topological type ${\mathbb R}^2$ may
flow to a flat surface under the hyperbolic geometric flow,
otherwise the flow may blow up in finite time and the surgery
technique has to be used. In some sense, we try to show the
following interesting phenomenon: for the hyperbolic geometric flow,
the surgery technique may be replaced by choosing suitable initial
velocity tensor.

The following \textit{Einstein's hyperbolic geometric flow} has also been
introduced by Kong and Liu \cite{kl}
\begin{equation}
\label{1.6}\dfrac{\partial^{2}g_{ij}}{\partial t^{2}}+2R_{ij}+\dfrac12g^{pq}
\dfrac{\partial g_{ij}}{\partial t}\dfrac{\partial g_{pq}}{\partial t}%
-g^{pq}\dfrac{\partial g_{ip}}{\partial t}\dfrac{\partial g_{qj}}{\partial
t}=0\qquad(i,j=1,2),
\end{equation}
where $g^{ij}$ is the inverse of $g_{ij}$. Noting the fact that the
equation (\ref{1.6}) is equivalent to (\ref{1.1}) for the metric
$g_{ij}$ with the form (\ref{1.5}), we know that all conclusions
mentioned above hold for the Einstein's hyperbolic geometric flow
(\ref{1.6}).

Here we would like to point out that, perhaps the method is more
important than the results obtained in this paper. Our method may
provide a new approach to the Penrose conjecture (see Penrose
\cite{p}) in general relativity and some of Yau's conjectures (e.g.,
problem 17 stated in Yau \cite{y}) about noncompact complete
manifolds with nonnegative curvature in differential geometry.

The paper is organized as follows. An interesting nonlinear partial
differential equation related to the metric (\ref{1.5}) is derived
in Section 2. Based on this, we prove the global existence theorem
on hyperbolic geometric flow, i.e., Theorem 1.1. The asymptotic
behavior theorem, i.e., Theorem 1.2 is proved in Section 3, the
proof depends on some new uniform \textit{a priori} estimates on
higher derivatives which are interesting in their own right. In
section 4, we investigate the blowup phenomena and the formation of
singularities in hyperbolic geometric flow. Section 5 is devoted to
the discussion about the radial solutions to the hyperbolic
geometric flow, i.e., the case $u=u(t,r)$ in which
$r=\sqrt{x^2+y^2}$. Section 6 concerns some geometric obstructions
to the existence of smooth long-time solutions and periodic
solutions of the hyperbolic geometric flow on general Riemann
surfaces.

\section{Global existence of hyperbolic geometric flow --- Proof of Theorem 1.1}

In our previous work \cite{dkl}, we have studied the flow of a metric by its
Ricci curvature
\[
\dfrac{\partial^{2}g_{ij}}{\partial t^{2}}=-2R_{ij}%
\]
on $n$-dimensional manifolds, where $n\geq5.$ In some respects the higher
dimensional cases are easier, due to the enough fast decay of the solution for
the corresponding linear wave equations. For surfaces, some estimates on the
curvature fails for that the solutions of the corresponding two-dimensional
linear wave equations only possesses ``slow" decay behavior. Therefore a new
approach is needed.

On a surface, the hyperbolic geometric flow equation simplifies, because all
of the information about curvature is contained in the scalar curvature
function $R.$ In our notation, $R=2K$ where $K$ is the Gauss curvature. The
Ricci curvature is given by
\begin{equation}
\label{2.1}R_{ij}=\dfrac12Rg_{ij},
\end{equation}
and the hyperbolic geometric flow equation simplifies the following equation
for the special metric
\begin{equation}
\label{2.2}\dfrac{\partial^{2}g_{ij}}{\partial t^{2}}=-Rg_{ij}.
\end{equation}
The metric for a surface can always be written (at least locally) in the
following form
\begin{equation}
\label{2.3}g_{ij}=u(t,x,y)\delta_{ij},
\end{equation}
where $u(t,x,y)>0.$ Therefore, we have
\begin{equation}
\label{2.4}R=-\dfrac{\triangle\ln u}{u}.
\end{equation}
Thus the equation (\ref{2.2}) becomes
\[
\dfrac{\partial^{2}u}{\partial t^{2}}=\dfrac{\triangle\ln u}{u}\cdot u,
\]
namely,
\begin{equation}
\label{2.5}u_{tt}-\triangle\ln u=0.
\end{equation}

In order to prove Theorem 1.1, by the uniqueness of the smooth solution of
nonlinear hyperbolic equations, it suffices to show that, for any given smooth
function $u_{1}(x)$ satisfying the conditions (a)-(b), the following Cauchy
problem
\begin{equation}
\label{2.61}\left\{
\begin{array}
[c]{l}%
\displaystyle{u_{tt}-\triangle\ln u=0,}\vspace{2mm}\\
\displaystyle{t=0:~~u=u_{0}(x,y),~u_{t}=u_{1}(x,y)}%
\end{array}
\right.
\end{equation}
has a unique solution for all time; Moreover, the derivatives of the solution
possess some algebraic decay estimates.

Notice that the initial data $u_{0}$ and $u_{1}$ only depend on the variable
$x$ and are independent of $y$. Therefore the Cauchy problem (\ref{2.61}) may
reduce to the following Cauchy problem for one-dimensional wave equation
\begin{equation}
\label{2.7}\left\{
\begin{array}
[c]{l}%
\displaystyle{u_{tt}- (\ln u)_{xx}=0,}\vspace{2mm}\\
\displaystyle{t=0:~~u=u_{0}(x),~u_{t}=u_{1}(x).}%
\end{array}
\right.
\end{equation}

In what follows, we shall solve the Cauchy problem (\ref{2.7}) and analyze its
solution's decay behavior.

Denote
\begin{equation}
\label{2.8}\phi=\ln u.
\end{equation}
Then (\ref{2.5}) reduces to
\begin{equation}
\label{2.9}\phi_{tt}-e^{-\phi}\triangle\phi=-\phi_{t}^{2}.
\end{equation}
In particular, the first equation in (\ref{2.7}) becomes%

\begin{equation}
\label{2.10}\phi_{tt}-e^{-\phi}\phi_{xx}=-\phi_{t}^{2}.
\end{equation}
Let
\begin{equation}
\label{2.11}v=\phi_{t},~~w=\phi_{x}.
\end{equation}
Then (\ref{2.10}) can be rewritten as the following quasilinear system of
first order
\begin{equation}
\label{2.12}\left\{
\begin{array}
[c]{l}%
\phi_{t}=v,\\
w_{t}-v_{x}=0,\\
v_{t}-e^{-\phi}w_{x}=-v^{2}%
\end{array}
\right.
\end{equation}
for smooth solutions.

Introduce
\begin{equation}
\label{2.13}p=v+e^{-\frac{\phi}{2}}w,~~q=v-e^{-\frac{\phi}{2}}w.
\end{equation}
We have

\begin{Lemma}
$p$ and $q$ satisfy
\begin{equation}
\label{2.14}\left\{
\begin{array}
[c]{c}%
\displaystyle{p_{t}-\lambda p_{x}=-\frac14\{p^{2}+3pq\},}\vspace{2mm}\\
\displaystyle{q_{t}+\lambda q_{x}=-\frac14\{q^{2}+3pq\},}%
\end{array}
\right.
\end{equation}
\end{Lemma}

where
\begin{equation}
\label{2.15}\lambda=e^{-\frac{\phi}{2}}. \quad\quad\quad\square
\end{equation}
\textbf{Proof.} We now calculate
\begin{equation}
\label{2.16}\begin{aligned}p_t-\lambda p_x&=\left(v+e^{-\frac{\phi}{2}}w\right)_t-e^{-\frac{\phi}{2}} \left(v+ e^{-\frac{\phi}{2}}w\right)_x\\ &=v_t+e^{-\frac{\phi}{2}}w_t-e^{-\frac{\phi}{2}}v_x-e^{-\phi}w_x-\frac12e^{-\frac{\phi}{2}}\phi_tw +\frac12e^{-\phi}\phi_xw\\ &=-v^2-\frac12vwe^{-\frac{\phi}{2}}+\frac12w^2e^{-\phi}\\ &=-v^2-\frac12v\left(we^{-\frac{\phi}{2}}\right)+\frac12\left(we^{-\frac{\phi}{2}}\right)^2\\ &=-\left(\frac{p+q}{2}\right)^2-\frac12\frac{p+q}{2}\frac{p-q}{2}+\frac12\left(\frac{p-q}{2}\right)^2\\ &=-\frac14\left(p^2+3pq\right).\end{aligned}
\end{equation}
In (\ref{2.16}) we have made use of (\ref{2.11}), (\ref{2.12}) and
(\ref{2.13}). In a similar way, we can prove the second equation in
(\ref{2.14}). Thus, the proof is finished.$\quad\quad\quad\blacksquare$

By a direct calculation, from (\ref{2.14}) we can obtain the following
interesting lemma.

\begin{Lemma}
It holds that
\begin{equation}
\label{2.17}\left\{
\begin{array}
[c]{l}%
\displaystyle{p_{t}-(\lambda p)_{x}=-pq,}\vspace{2mm}\\
\displaystyle{q_{t}+(\lambda q)_{x}=-pq.}%
\end{array}
\right.  \quad\quad\quad\square
\end{equation}

\end{Lemma}

Noting (\ref{2.7})-(\ref{2.8}), (\ref{2.11}) and (\ref{2.13}), we denote
\begin{equation}
\label{2.18}p_{0}(x)\triangleq\frac{u_{1}(x)}{u_{0}(x)}+\frac{u_{0}^{\prime
}(x)}{u_{0}^{\frac32}(x)},~~~ q_{0}(x)\triangleq\frac{u_{1}(x)}{u_{0}%
(x)}-\frac{u_{0}^{\prime}(x)}{u_{0}^{\frac32}(x)}.
\end{equation}
It is easy to verify that the $C^{2}$ solution of the Cauchy problem
(\ref{2.7}) is equivalent to the $C^{1}$ solution of the Cauchy problem for
the following quasilinear system of first order
\begin{equation}
\label{2.19}\left\{
\begin{array}
[c]{l}%
\displaystyle{\phi_{t}=\frac{p+q}{2},}\vspace{2mm}\\
\displaystyle{p_{t}-\lambda p_{x}=-\frac14(p^{2}+3pq),}\vspace{2mm}\\
\displaystyle{q_{t}+\lambda q_{x}=-\frac14(q^{2}+3pq)}%
\end{array}
\right.
\end{equation}
with the initial data
\begin{equation}
\label{2.20}t=0:\quad\phi=\ln u_{0}(x),~p=p_{0}(x),~q=q_{0}(x),
\end{equation}
where $\lambda=\lambda(\phi)$ is defined by (\ref{2.15}), $p_{0}(x)$ and
$q_{0}(x)$ are given by (\ref{2.18}). Obviously, (\ref{2.19}) is a strictly
hyperbolic system with three distinguished characteristics
\begin{equation}
\label{2.21}\lambda_{1}=-\lambda,\quad \lambda_{2}=0,\quad
\lambda_{3}=\lambda.
\end{equation}

In order to prove Theorem 1.1, it suffices to show the following
theorem.

\begin{Theorem}
If $u_{1}(x)$ is a smooth function with bounded $C^{1}$ norm and satisfies
\begin{equation}
\label{2.22}u_{1}(x)\geqslant\frac{|u_{0}^{\prime}(x)|}{\sqrt{u_{0}(x)}%
},~~\forall~x\in\mathbb{R},
\end{equation}
then the Cauchy problem (\ref{2.19}), (\ref{2.20}) has a unique global smooth
solution for all time $t\in\mathbb{R}$.$\quad\quad\quad\square$
\end{Theorem}

\begin{Corollary}
Under the assumptions in Theorem 2.1, the Cauchy problem (\ref{2.7}) has a
unique global $C^{2}$ solution for all time $t\in\mathbb{R}.\quad\quad
\quad\square$
\end{Corollary}

According to the existence and uniqueness theorem of smooth solution of
hyperbolic systems of first order, there exists a locally smooth solution of
the Cauchy problem (\ref{2.19})--(\ref{2.20}). In order to prove Theorem 2.1,
it suffices to establish the uniform \textit{a priori} estimate on the $C^{1}$
norm of $(\phi,p,q)$ in the domain where the smooth solution of the Cauchy
problem (\ref{2.19})--(\ref{2.20}) exists. That is to say, we have to
establish the uniform \textit{a priori} estimates on the $C^{0}$ norm of
$(\phi,p,q)$ and their derivatives of first order on the existence domain of
smooth solution of the Cauchy problem (\ref{2.19})--(\ref{2.20}). Noting
(\ref{2.19}), we see that the key point is to establish \textit{a priori}
estimate on the $C^{1}$ norm of $p$ and $q$.

In order to prove Theorem 2.1, we need the following lemmas.

\begin{Lemma}
In the existence domain of the smooth solution of the Cauchy problem
(\ref{2.19})--(\ref{2.20}), it holds that
\begin{equation}
\label{2.23}0\leqslant p(t,x)\leqslant\sup_{y\in\mathbb{R}} p_{0}(y)
\end{equation}
and
\begin{equation}
\label{2.24}0\leqslant q(t,x)\leqslant\sup_{y\in\mathbb{R}} q_{0}%
(y).\quad\quad\quad\square
\end{equation}

\end{Lemma}

\noindent\textbf{{Proof.}} In fact, passing through any point $(t,x)$, we can
draw two characteristics, defined by $\xi=\xi_{\pm}(\tau;t,x)$, which satisfy
\begin{equation}
\label{2.25}\left\{
\begin{array}
[c]{l}%
\frac{d\xi_{\pm}}{d\tau}=\pm\lambda(\tau;\xi_{\pm}(\tau;t,x)), \vspace{2mm}\\
\xi_{\pm}(t;t,x)=x,
\end{array}
\right.
\end{equation}
respectively. Noting the last two equations in (\ref{2.19}), we observe that,
along the characteristic $\xi=\xi_{+}(\tau;t,x)$, it holds that
\begin{equation}
\label{2.26}p(t,x)=p_{0}(\xi_{+}(0;t,x))\exp\left\{  \int_{0}^{t}%
-\frac14\left[  p+3q\right]  (\tau;\xi_{+} (\tau;t,x))d\tau\right\}  .
\end{equation}

On the other hand, noting (\ref{2.18}) and (\ref{2.20}), we have
\begin{equation}
\label{2.27}p_{0}(x)\geqslant0,~~\forall\;x\in\mathbb{R},
\end{equation}
and then by (\ref{2.26}), we obtain
\[
p(t,x)\geqslant0
\]
in the existence domain of the smooth solution of the Cauchy problem
(\ref{2.19})-(\ref{2.20}).

Similarly, we can prove
\begin{equation}
\label{2.28}q(t,x)\geqslant0.
\end{equation}

On the other hand, noting (\ref{2.27}), we obtain from (\ref{2.26}) that, for
any point $(t,x)$ in the existence domain of the smooth solution
\begin{equation}
\label{2.29}0\leqslant p(t,x)\leqslant p_{0}(\xi_{+}(0;t,x))\leqslant
\sup_{y\in\mathbb{R}}p_{0}(y).
\end{equation}
This is the desired inequality (\ref{2.23}).

Similarly, we can prove (\ref{2.24}). Thus, the proof is completed.$\quad
\quad\quad\blacksquare$

\vskip 4mm

We next estimate $p_{x}$ and $q_{x}$.

Let
\begin{equation}
\label{2.30}r=p_{x},~~s=q_{x}.
\end{equation}
Similar to Lemma 2.1, we have

\begin{Lemma}
$r$ and $s$ satisfy
\begin{equation}
\label{2.31}\left\{
\begin{array}
[c]{l}%
r_{t}-\lambda r_{x}=-\frac14\left[  (2q+3p)r+3ps\right]  ,\vspace{2mm}\\
s_{t}+\lambda s_{x}=-\frac14\left[  (2p+3q)s+3qr\right]  .
\end{array}
\right.  \quad\quad\quad\square
\end{equation}
\end{Lemma}

\noindent\textbf{Proof.} By a direct calculation, we can easily prove
(\ref{2.31}). $\quad\quad\quad\blacksquare$

\vskip 4mm

Denote
\begin{equation}
\label{2.32}r_{0}(x)\triangleq p_{0}^{\prime}(x),~~~~s_{0}(x)\triangleq
q_{0}^{\prime}(x).
\end{equation}
We have

\begin{Lemma}
In the existence domain of the smooth solution, it holds that
\begin{equation}
\label{2.33}|r(t,x)|,~~|s(t,x)|\leqslant\max\left\{  \sup_{y\in\mathbb{R}%
}|r_{0}(y)|,~~ \sup_{y\in\mathbb{R}}|s_{0}(y)|\right\}  .\quad\quad
\quad\square
\end{equation}
\end{Lemma}

\noindent\textbf{Proof.} Let
\begin{equation}
\label{2.34}A=\frac{2q+3p}{4},\quad B=\frac{3p}{4}, \quad\bar{A}=\frac
{2p+3q}{4},\quad\bar{B}=\frac{3q}{4}.
\end{equation}
Then the system (\ref{2.31}) can be rewritten as
\begin{equation}
\label{2.35}\left\{
\begin{array}
[c]{l}%
r_{t}-\lambda r_{x}=-Ar-Bs,\vspace{2mm}\\
s_{t}+\lambda s_{x}=-\bar{A}s-\bar{B}r.
\end{array}
\right.
\end{equation}
By Lemma 2.3, we have
\begin{equation}
\label{2.36}A,\;\bar{A},\;B,\;\bar{B}\geqslant0,\quad A\geqslant B,
\quad\mathrm{and}\quad\bar{A}\geqslant\bar{B}.
\end{equation}
By the terminology in Kong \cite{k}, the system (\ref{2.35}) (i.e.,
(\ref{2.31})) is \textit{weakly dissipative}. Therefore, it follows from
Theorem 2.3 in Kong \cite{k} that
\begin{equation}
\label{2.37}|r(t,x)|,~~|s(t,x)|\leqslant\max\left\{  \sup_{y\in I(t,x)}%
|r_{0}(y)|,~~ \sup_{y\in I(t,x)}|s_{0}(y)|\right\}  ,
\end{equation}
where
\[
I(t,x)=\left[  \xi_{-}(0;t,x),\xi_{+}(0;t,x)\right]  .
\]
The desired estimate (\ref{2.33}) comes from (\ref{2.37}) directly. This
proves Lemma 2.5. $\quad\quad\quad\blacksquare$

\vskip 4mm

We next estimate $\phi$.

For any fixed point $(t,x)$ in the existence domain of the smooth solution, it
follows from the first equation in (\ref{2.19}) that
\begin{equation}
\label{2.38}\phi(t,x)=\phi(0,x)+\int_{0}^{t}\frac{p+q}{2}(\tau,x)d\tau.
\end{equation}
Noting (\ref{2.20}) and (\ref{2.23})-(\ref{2.24}), we have
\begin{equation}
\label{2.39}\ln u_{0}(x)\leqslant\phi(t,x)\leqslant\ln u_{0}(x)+\frac
12\left\{  \sup_{y\in\mathbb{R}}|p_{0}(y)|+\sup_{y\in\mathbb{R}}%
|q_{0}(y)|\right\}  t.
\end{equation}
On the other hand, deriving the first equation in (\ref{2.19}) with respect to
$x$ gives
\begin{equation}
\label{2.40}(\phi_{x})_{t}=\frac12(p_{x}+q_{x})=\frac12(r+s),
\end{equation}
that is,
\begin{equation}
\label{2.41}\phi_{x}=\phi_{x}(0,x)+\frac12\int_{0}^{t}(r+s)(\tau,x)d\tau.
\end{equation}
Using (\ref{2.20}) and noting (\ref{2.37}), we have
\begin{equation}
\label{2.42}|\phi_{x}(t,x)|\leqslant\frac{|u_{0}^{\prime}(x)|}{u_{0}(x)}+
\left[  \sup_{y\in\mathbb{R}}\left|  r_{0}(y)\right|  +\sup_{y\in\mathbb{R}%
}|s_{0}(y)|\right]  t.
\end{equation}

\noindent\textbf{Proof of Theorem 2.1.} Notice that $u_{0}(x)$ is a smooth
function with bounded $C^{2}$ norm and satisfies (\ref{1.3}). On the other
hand, notice that $u_{1}(x)$ is a smooth function with bounded $C^{1}$ norm
and satisfies (\ref{2.22}). Then it follows from (\ref{2.18}) that
\begin{equation}
\label{2.43}\left\{
\begin{array}
[c]{l}%
\displaystyle{0\leqslant p_{0}(x)\leqslant\sup_{x\in\mathbb{R}} p_{0}%
(x)\triangleq P_{0}<\infty,}\vspace{2mm}\\
\displaystyle{0\leqslant q_{0}(x)\leqslant\sup_{x\in\mathbb{R}} q_{0}%
(x)\triangleq Q_{0}<\infty.}%
\end{array}
\right.
\end{equation}
Thus, it follows from (\ref{2.39}) that, in the existence of the smooth
solution
\begin{equation}
\label{2.44}\ln m\leqslant\phi(t,x)\leqslant\ln M+\frac12(P_{0}+Q_{0})t.
\end{equation}
On the other hand, noting (\ref{2.15}), we have
\begin{equation}
\label{2.45}M^{-\frac12}e^{-\frac14(P_{0}+Q_{0})t}\leqslant\lambda\leqslant
m^{-\frac12}.
\end{equation}

>From any interval $[a,b]$ in the $x$-axis, we introduce the
following triangle domain, which is the strong determinate domain of
the interval $[a,b]$
\begin{equation}
\label{2.46}\triangle_{[a,b]}=\left\{  (t,x)\left|  a+ m^{-\frac12}t\leqslant
x\leqslant b+M^{-\frac12}t\right.  \right\}  .
\end{equation}

It is easy to see that, in order to prove Theorem 2.1, it suffices to prove
that, for arbitrary interval $[a,b]$ in the $x$-axis the Cauchy problem
(\ref{2.19})-(\ref{2.20}) has a unique smooth solution on $\triangle_{[a,b]}$.

Noting the estimates (\ref{2.23})-(\ref{2.24}), (\ref{2.37}), (\ref{2.39}) and
(\ref{2.42}), for any point $(t,x)$ in $\triangle_{[a,b]}$ we have following a
priori estimates
\begin{equation}
\label{2.47}0\leqslant p(t,x)\leqslant\sup_{y\in[a,b]}p_{0}(y)\le
P_{0},
\end{equation}
\begin{equation}
\label{2.48}0\leqslant q(t,x)\leqslant\sup_{y\in[a,b]}q_{0}(y)\le
Q_{0},
\end{equation}
\begin{equation}
\label{2.49}|r(t,x)|,~|s(t,x)|\leqslant\max\left\{  \sup_{y\in[a,b]}%
|r_{0}(y)|, \sup_{y\in[a,b]}|s_{0}(y)|\right\}  <\infty,
\end{equation}
\begin{equation}
\label{2.50}\ln m\leqslant\phi(t,x)\leqslant\ln M+\frac12\left\{  \sup
_{y\in[a,b]}|p_{0}(y)|+\sup_{y\in[a,b]}|q_{0}(y)| \right\}  t\leqslant\ln
M+\frac12(P_{0}+Q_{0})t_{[a,b]},
\end{equation}
\begin{equation}
\label{2.51}|\phi_{x}(t,x)|\leqslant\frac{\displaystyle\sup_{x\in[a,b]}%
|u_{0}^{\prime}(x)|}{m}+ \left\{  \sup_{y\in[a,b]}|r_{0}(y)|+\sup_{y\in
[a,b]}|s_{0}(y)|\right\}  t_{[a,b]},
\end{equation}
where
\[
t_{[a,b]}=\frac{b-a}{m^{-\frac12}-M^{-\frac12}}.
\]
In the estimates (\ref{2.47})-(\ref{2.51}), we assume that the solution
exists. The above \textit{a priori} estimates (\ref{2.47})-(\ref{2.51})
implies that the Cauchy problem (\ref{2.19})-(\ref{2.20}) has a unique smooth
solution on the whole triangle domain $\triangle_{[a,b]}$. This proves Theorem
2.1. $\quad\quad\quad\blacksquare$

\vskip 4mm

Theorem 1.1 follows from Theorem 2.1 immediately.

\section{Asymptotic behavior --- Proof of Theorem 1.2}

In this section, we shall establish some uniform estimates on the
global smooth solution $(\phi,p,q)$ of the Cauchy problem
(\ref{2.19})-(\ref{2.20}) as well as its derivatives $r$ and $s$.
Based on this, we can prove Theorem 1.2.

We first prove the part (I) of Theorem 1.2.

We now assume that (1.4) holds, i.e.,
\begin{equation}
\label{3.1}u_{1}(x)\geqslant\frac{|u_{0}^{\prime}(x)|}{\sqrt{u_{0}(x)}%
},~~~\forall\; x\in\mathbb{R}.
\end{equation}
Therefore, all the estimates mentioned in the previous section hold.
Under the hypothesis (\ref{3.1}), we next establish some decay
estimates which play an important role in the proof of Theorem 1.2.

Noting the assumption (\ref{3.1}) (i.e., (\ref{1.3})), we obtain
from (\ref{2.18}) that
\begin{equation}
\label{3.2}p_{0}(x),\;\;q_{0}(x)\geqslant 0,~~\forall\;
x\in\mathbb{R}.
\end{equation}
In the present situation, similar to the argument in Section 2, we have
\begin{equation}
\label{3.3}p(t,x)\ge 0,\;\;q(t,x)\ge 0,\quad\forall\; (t,x)\in\mathbb{R}^{+}%
\times\mathbb{R}.
\end{equation}
Thus, it follows from (\ref{2.19}) that
\begin{equation}
\label{3.4}p_{t}-\lambda p_{x}\leqslant-\frac14p^{2}%
\end{equation}
and
\begin{equation}
\label{3.5}q_{t}+\lambda q_{x}\leqslant-\frac14q^{2}.
\end{equation}
For any fixed point $(0,\alpha)$ in the $x$-axis, along the characteristic
$\xi=\xi_{-}(t;0,\alpha)$ it holds that
\begin{equation}
\label{3.6}p(t,\xi_{-}(t;0,\alpha))\leqslant\dfrac{p_0(\alpha)}{1+\dfrac14p_{0}%
(\alpha)t}\leqslant P_0,\quad\forall\; t\geqslant0.
\end{equation}
This implies that, for any point $(t,x)\in\mathbb{R}^{+}\times\mathbb{R}$, we
have
\begin{equation}
\label{3.7}p(t,x)\leqslant\dfrac{p_{0}(\xi_{-}(0;t,x))}{1+\dfrac14p_{0}(\xi_{-}(0;t,x))t}%
\leqslant P_0,
\end{equation}
namely,
\begin{equation}
\label{3.8} p(t,x)\leqslant C_{1}, \quad\forall\; (t,x)\in
\mathbb{R}^{+}\times\mathbb{R},\end{equation} here and hereafter
$C_{i}\; (i=1,2,\cdots)$ stand for the positive constants
independent of $t$ and $x$, but depending on $\varepsilon, M$ and
the $C^{2}$ norm of $u_{0}$ and $C^{1}$ norm of $u_{1}$.

Similarly, we can prove
\begin{equation}
\label{3.9} q(t,x)\leqslant C_2, \quad\forall\;
(t,x)\in\mathbb{R}^{+}\times\mathbb{R}.\end{equation}

In what follows, by maximum principle for hyperbolic systems (see
Theorem 2.1 in Kong \cite{k}), we establish the decay estimates on
$r$ and $s$.

Noting (\ref{2.34}), we have
\begin{equation}
\label{3.10}A-B=\frac12q\ge 0.
\end{equation}
We now apply the maximum principle derived in Kong \cite{k} to the system
(\ref{2.35}) and obtain
\begin{equation}
\label{3.11}|r(t,x)|\leqslant\max\left\{  \sup_{y\in
I(t,x)}|r_{0}(y)|, \sup_{y\in I(t,x)}|s_{0}(y)|\right\},
\end{equation}
where $I(t,x)$ and $\xi_{-}(\tau;t,x)$ are defined as before. Let
\begin{equation}
\label{3.12}N=\sup_{y\in\mathbb{R}}|r_{0}(y)|+\sup_{y\in\mathbb{R}}|s_{0}(y)|.
\end{equation}
By (\ref{3.3}), it follows from and (\ref{3.11}) that
\begin{equation}
\label{3.13}\begin{aligned}|r(t,x)|\leqslant N,
\quad\forall\;(t,x)\in\mathbb{R}^+\times\mathbb{R}.\end{aligned}
\end{equation}

Similarly, we have
\begin{equation}
\label{3.14}|s(t,x)|\leqslant N,
\quad\forall\;(t,x)\in\mathbb{R}^+\times\mathbb{R}.
\end{equation}

Next we estimate $\phi(t,x)$ and its derivatives. It follows from
(\ref{2.39}) that
\begin{equation}
\label{3.15}\ln m\leqslant\phi(t,x)\leqslant\ln M+\frac12(P_{0}+Q_{0})t,
\end{equation}
as a consequence, by (\ref{2.8})
\begin{equation}
\label{3.16}m\leqslant u(t,x)\leqslant M\exp\left\{  \frac12(P_{0}%
+Q_{0})t\right\}  .
\end{equation}
To estimate $\phi_{t}$ and $\phi_{x}$, we use (\ref{2.11}) and (\ref{2.13}).
By (\ref{2.11}) and (\ref{2.13}), we have
\begin{equation}
\label{3.17}\phi_{t}=\frac12(p+q),\quad\phi_{x}=\frac12e^{\frac{\phi}{2}%
}(p-q).
\end{equation}
Thus, it follows from (\ref{3.8}) and (\ref{3.9}) that
\begin{equation}
\label{3.18}|\phi_{t}(t,x)|\leqslant\frac{C_{1}+C_{2}}{2}\leqslant
C_{3},
\end{equation}
\begin{equation}
\label{3.19}|\phi_{x}(t,x)|\leqslant\frac12\exp\left\{ \frac12\left[
\ln M+\frac12(P_{0}+Q_{0})t\right]  \right\}
\frac{C_{1}+C_{2}}{2}\leqslant C_{4}\exp\{C_{5}t\}.
\end{equation}

We now estimate $\phi_{xx}$. Noting (\ref{2.30}) and (\ref{2.11}), we obtain
from (\ref{2.13}) that
\begin{equation}
\label{3.20}\left\{
\begin{array}
[c]{l}%
\displaystyle{r=p_{x}=v_{x}+e^{-\frac{\phi}{2}}\phi_{xx}-\frac12
e^{-\frac{\phi}{2}}\phi_{x}^{2},}\vspace{2mm}\\
\displaystyle{s=q_{x}=v_{x}-e^{-\frac{\phi}{2}}\phi_{xx}+\frac12
e^{-\frac{\phi}{2}}\phi_{x}^{2}.}%
\end{array}
\right.
\end{equation}
It follows from (\ref{3.20}) that
\begin{equation}
\label{3.21}\phi_{xx}=\frac12e^{\frac{\phi}{2}}\left\{  (r-s)+e^{-\frac{\phi
}{2}}\phi_{x}^{2}\right\}  =\frac12\left\{  (r-s)e^{\frac{\phi}{2}}+\phi
_{x}^{2}\right\}  .
\end{equation}
Thus, by (\ref{3.15}), (\ref{3.13})--(\ref{3.14}) and (\ref{3.19}) we obtain
\begin{equation}
\label{3.22}|\phi_{xx}(t,x)|\leqslant C_{6}\exp\left\{
C_{7}t\right\} .
\end{equation}

We finally estimate the scalar curvature $R$. By the definition, we
have
\begin{equation}
\label{3.23}R=\frac{(\ln u)_{xx}}{u}.
\end{equation}
Noting (\ref{2.8}), we obtain from (\ref{3.23}) that
\begin{equation}
\label{3.24}R(t,x)=\frac{\phi_{xx}}{e^{\phi}}.
\end{equation}
Using (\ref{3.21}), we have
\begin{equation}
\label{3.25}R(t,x)=\frac{\frac12\left\{  (r-s)e^{\frac{\phi}{2}}+\phi_{x}%
^{2}\right\}  } {e^{\phi}}=\frac12\left\{  (r-s)e^{-\frac{\phi}{2}}+\left(
\phi_{x}e^{-\frac{\phi}{2}}\right)  ^{2}\right\}  .
\end{equation}
Thus, by (\ref{2.13}) we get
\begin{equation}
\label{3.26}R(t,x)=\frac12\left\{  (r-s)e^{-\frac{\phi}{2}}+\left(  \frac
{p-q}{2}\right)  ^{2}\right\}  .
\end{equation}
Then, using (\ref{3.13})--(\ref{3.15}) and (\ref{3.8})--(\ref{3.9}), we have
\begin{equation}
\label{3.27}\begin{aligned}|R(t,x)|&\leqslant\frac12\left\{2N\cdot
\exp\left\{-\frac12\ln m\right\}+\frac14(C_1+C_2)^2\right\}
\leqslant C_{8}, \quad\forall\;
(t,x)\in\mathbb{R}^+\times\mathbb{R}.\end{aligned}
\end{equation}
(\ref{3.27}) implies that, for any fixed
$(t,x)\in\mathbb{R}^2\times\mathbb{R}$ there exists a positive
constant $k$ such that
\begin{equation}
\label{3.28}|R(t,x)|\le k,
\end{equation}
where $k$ depends on $m,\;M$, the $C^{2}$ norm of $u_{0}$ and
$C^{1}$ norm of $u_{1}$, but independent of $t$ and $x$. This proves
the part (I) of Theorem 1.2.

\begin{Remark}
>From the above argument we observe that, in order to estimate the
scalar curvature $R$, we need to estimate the quantities:
$re^{-\frac{\phi}{2}}$ and $se^{-\frac{\phi}{2}}$. A more convenient
way is as follows: let
\begin{equation}\tilde{r}=re^{-\frac{\phi}{2}}, \quad
\tilde{s}=-se^{-\frac{\phi}{2}},\end{equation} then it follows from
(\ref{2.31}) that
\begin{equation}
\left\{
\begin{array}
[c]{l}%
\tilde{r}_{t}-\lambda \tilde{r}_{x}=-\frac34 p(\tilde{r}-\tilde{s}) ,\vspace{2mm}\\
\tilde{s}_{t}+\lambda \tilde{s}_{x}=-\frac34 q(\tilde{s}-\tilde{r}).
\end{array}
\right.
\end{equation}
Noting (3.3) and using Remark 4 in Hong \cite{h} or Theorem 2.4 in
Kong \cite{k}, we have
\begin{equation}
\min\{\inf_{x\in\mathbb{R}}\tilde{r}(0,x),\inf_{x\in\mathbb{R}}\tilde{s}(0,x)\}\le
\tilde{r}(t,x),\;\tilde{s}(t,x)\le
\max\{\sup_{x\in\mathbb{R}}\tilde{r}(0,x),\sup_{x\in\mathbb{R}}\tilde{s}(0,x)\}
\end{equation}
for any $(t,x)\in \mathbb{R}^+\times \mathbb{R}$. (3.31) gives the
uniform estimates on the upper and lower bounds of the quantities
$re^{-\frac{\phi}{2}}$ and $se^{-\frac{\phi}{2}}$.
\end{Remark}

We next prove the part (I\!I) of Theorem 1.2.

In what follows, we assume that (1.4a) is satisfied, i.e., it holds
that
$$u_{1}(x)\geqslant\frac{|u_{0}^{\prime}(x)|}{\sqrt{u_{0}(x)}%
}+\varepsilon,~~~\forall\; x\in\mathbb{R},\eqno{(1.4a)}$$ where
$\varepsilon$ is an arbitrary positive constant. Obviously, under
the assumption (1.4a), (\ref{2.22}) is always true, and then all the
estimates mentioned above hold. In fact, under the hypothesis
(1.4a), we may furthermore establish some decay estimates which play
an important role in the proof of the part (I\!I) of Theorem 1.2.

Noting the assumption (1.4a), we obtain from (\ref{2.18}) that
\begin{equation}
\label{3.32}p_{0}(x),\;\;q_{0}(x)\geqslant\frac{\varepsilon}{M}\triangleq
M^{\prime},~~\forall\; x\in\mathbb{R}.
\end{equation}

Introduce
\begin{equation}
\label{3.33}\hat{p}=up,\quad \hat{q}=uq.
\end{equation}
Then it follows from (2.8), (2.11) and (2.13) that
\begin{equation}
\label{3.34} \hat{p}=u_t+\frac{1}{\sqrt{u}}u_x,\quad
\hat{q}=u_t-\frac{1}{\sqrt{u}}u_x.
\end{equation}
This gives
\begin{equation}
\label{3.35} u(t,x)=u_0(x)+\int_0^t(\hat{p}+\hat{q})(\tau,x)d\tau.
\end{equation}
On the other hand, similar to (2.14), we have
\begin{equation}
\label{3.36}\left\{
\begin{array}
[c]{c}%
\displaystyle{\hat{p}_{t}- \frac{1}{\sqrt{u}}\hat{p}_{x}=\frac{1}{4u}\hat{p}(\hat{q}-\hat{p}),}\vspace{2mm}\\
\displaystyle{\hat{q}_{t}+ \frac{1}{\sqrt{u}}\hat{q}_{x}=\frac{1}{4u}\hat{q}(\hat{p}-\hat{q}).}%
\end{array}
\right.
\end{equation}
Noting (3.3), (3.16) and (3.33), we observe that
\begin{equation}
\label{3.37}\frac{1}{4u}\hat{p}\ge 0\quad \frac{1}{4u}\hat{q}\ge 0,
\quad\forall\; (t,x)\in\mathbb{R}^+\times\mathbb{R}.
\end{equation}
Then, by Remark 4 in Hong \cite{h} or Theorem 2.4 in Kong \cite{k},
it follows from (3.36) that
\begin{equation}\label{3.38}
\min\{\inf_{x\in\mathbb{R}}\hat{q}(0,x),\inf_{x\in\mathbb{R}}\hat{q}(0,x)\}\le
\hat{p}(t,x),\;\hat{q}(t,x)\le
\max\{\sup_{x\in\mathbb{R}}\hat{p}(0,x),\sup_{x\in\mathbb{R}}\hat{q}(0,x)\}
\end{equation}
for any $(t,x)\in \mathbb{R}^+\times \mathbb{R}$. Noting (1.3) and
(3.32), we have
\begin{equation}
\label{3.39}C_9\le \hat{p}(t,x),\;\; \hat{q}(t,x)\le C_{10},
\quad\forall\; (t,x)\in\mathbb{R}^+\times\mathbb{R}.
\end{equation}
Thus, it follows from (3.35) that
\begin{equation}
\label{3.40}C_{11}(1+t)\le u(t,x)\le C_{12}(1+t), \quad\forall\;
(t,x)\in\mathbb{R}^+\times\mathbb{R}.
\end{equation}
Noting (3.33) gives
\begin{equation}
\label{3.41}\frac{C_{13}}{1+t}\le p(t,x),\;\; q(t,x)\le
\frac{C_{14}}{1+t}, \quad\forall\;
(t,x)\in\mathbb{R}^+\times\mathbb{R},
\end{equation}
where $C_{13}\le C_9/C_{12}$ and $C_{14}\ge C_{10}/C_{11}$

\vskip 4mm

We next establish some decay estimates on $r$ and $s$.

Noting (\ref{2.34}), we have
\begin{equation}
\label{3.42}A-B=\frac12q.
\end{equation}
We now apply the maximum principle derived in Kong \cite{k} to the
system (\ref{2.35}) and obtain
\begin{equation}
\label{3.43}|r(t,x)|\leqslant\max\left\{  \sup_{y\in
I(t,x)}|r_{0}(y)|, \sup_{y\in I(t,x)}|s_{0}(y)|\right\}  \exp\left\{
-\int_{0}^{t}\dfrac 12q(\tau,\xi_{-}(\tau;t,x))d\tau\right\}  ,
\end{equation}
where $I(t,x)$ and $\xi_{-}(\tau;t,x)$ are defined as in Section 2.

Let
\begin{equation}
\label{3.44}N=\sup_{y\in\mathbb{R}}|r_{0}(y)|+\sup_{y\in\mathbb{R}}|s_{0}(y)|.
\end{equation}
Noting (3.41), we obtain from (3.43) that
\begin{equation}
\label{3.45}\begin{aligned}|r(t,x)|&\leqslant
N\exp\left\{-\frac{C_{13}}{2}\int_0^t \frac{1}{1+\tau}d\tau\right\}\\
&\leqslant
N\exp\left\{\ln(1+t)^{-C_{15}}\right\}\\&=\frac{N}{(1+t)^{C_{15}}}.\end{aligned}
\end{equation}
Similarly, we have
\begin{equation}
\label{3.46}|s(t,x)|\leqslant\frac{N}{(1+t)^{C_{16}}}.
\end{equation}

We now estimate $\phi_{xx}$.

Noting (\ref{2.30}) and (\ref{2.11}), we obtain from (\ref{2.13})
that
\begin{equation}
\label{3.47}\left\{
\begin{array}
[c]{l}%
\displaystyle{r=p_{x}=v_{x}+e^{-\frac{\phi}{2}}\phi_{xx}-\frac12
e^{-\frac{\phi}{2}}\phi_{x}^{2},}\vspace{2mm}\\
\displaystyle{s=q_{x}=v_{x}-e^{-\frac{\phi}{2}}\phi_{xx}+\frac12
e^{-\frac{\phi}{2}}\phi_{x}^{2}.}%
\end{array}
\right.
\end{equation}
It follows from (\ref{3.47}) that
\begin{equation}
\label{3.48}\phi_{xx}=\frac12e^{\frac{\phi}{2}}\left\{
(r-s)+e^{-\frac{\phi }{2}}\phi_{x}^{2}\right\}  =\frac12\left\{
(r-s)e^{\frac{\phi}{2}}+\phi _{x}^{2}\right\}  .
\end{equation}

We finally estimate the scalar curvature $R$.

By the definition, we have
\begin{equation}
\label{3.49}R=\frac{(\ln u)_{xx}}{u}.
\end{equation}
Noting (\ref{2.8}), we obtain from (\ref{3.49}) that
\begin{equation}
\label{3.50}R(t,x)=\frac{\phi_{xx}}{e^{\phi}}.
\end{equation}
Using (\ref{3.48}), we have
\begin{equation}
\label{3.51}R(t,x)=\frac{\frac12\left\{  (r-s)e^{\frac{\phi}{2}}+\phi_{x}%
^{2}\right\}  } {e^{\phi}}=\frac12\left\{
(r-s)e^{-\frac{\phi}{2}}+\left( \phi_{x}e^{-\frac{\phi}{2}}\right)
^{2}\right\}  .
\end{equation}
Noting (2.8), (2.11) and (\ref{2.13}), we get
\begin{equation}
\label{3.52}R(t,x)=\frac12\left\{\frac{r-s}{\sqrt{u}}+\left( \frac
{p-q}{2}\right)  ^{2}\right\}  .
\end{equation}
Then, using (\ref{3.40})-(\ref{3.41}) and (\ref{3.45})-(\ref{3.46}),
we have
\begin{equation}
\label{3.53}\begin{aligned}|R(t,x)|&\leqslant\frac12\left\{\frac{2N}{(1+t)
^{\min\{C_{15},C_{16}\}}\sqrt{C_{11}(1+t)}}+\left(\frac{C_{14}-C_{13}}{2(1+t)}\right)^2\right\}
\\&\leqslant C_{17}\left\{\frac{1}{(1+t)^{\min\{C_{15},C_{16}\}+\frac12}}+\frac{1}{(1+t)^2}\right\}\\ &
\leqslant
C_{18}\frac{1}{(1+t)^{\min\{2,\min\{C_{15},C_{16}\}+\frac12\}}},~~~~~\forall\;
(t,x)\in\mathbb{R}^+\times\mathbb{R}.\end{aligned}
\end{equation}
Taking $\tilde{k}\ge C_{18}$ and
$\gamma=\min\{2,\min\{C_{15},C_{16}\}+\frac12\}$, we obtain from
(3.53) that
\begin{equation}
\label{3.54}|R(t,x)|\le
\dfrac{\tilde{k}}{(1+t)^{\gamma}},\quad\forall\;
(t,x)\in\mathbb{R}^+\times\mathbb{R}.
\end{equation}
This proves the part (I\!I) of Theorem 1.2. Thus, the proof of
Theorem 1.2 is completed.

\section{Blowup phenomena and formation of singularities --- Proof of Theorem 1.3}

In this section we will investigate the blowup phenomena of hyperbolic
geometric flow and the formation of singularities, provided that the
assumption (\ref{1.3}) is not satisfied.

Throughout this section, we assume that (\ref{1.2}) holds. As in
Theorem 1.3, we assume that there exists a point
$x_{0}\in\mathbb{R}$ such that
\begin{equation}
\label{4.1}u_{0}^{\prime}(x_{0})<0.
\end{equation}
In this case, we choose
\begin{equation}
\label{4.2}u_{1}(x)\equiv\frac{u_{0}^{\prime}(x)}{\sqrt{u_{0}(x)}}%
,~~~\forall\; x\in\mathbb{R}.
\end{equation}
In what follows, we shall prove

\begin{Theorem}
For the initial data $u_{0}(x)$ and $u_{1}(x)$ mentioned above, the Cauchy
problem (\ref{1.4}) has a unique smooth solution only in $[0,\tilde{T}_{\max
})\times\mathbb{R}$, where
\begin{equation}
\label{4.3}\tilde{T}_{\max}=-\frac{4}{\inf_{x\in\mathbb{R}}\left\{
p_{0}(x)\right\}  },
\end{equation}
where
\[
p_{0}(x)=2u_{0}^{\prime}(x)u_{0}^{-\frac{3}{2}}(x).
\]
Moreover, there exists some point $(\tilde{T}_{\max},x_{\ast})$ such that the
scalar curvature $R(t,x)$ satisfies
\begin{equation}
\label{4.4}R(t,x)\rightarrow\infty\qquad\text{as}~~(t,x)\nearrow(\tilde
{T}_{\max},x_{\ast}). \quad\quad\quad\square
\end{equation}

\end{Theorem}

\noindent\textbf{Proof.} As in Sections 2-3, it suffices to study the Cauchy
problem (\ref{2.7}), equivalently, the Cauchy problem
\begin{equation}
\label{4.5}\left\{
\begin{array}
[c]{l}%
\displaystyle{\phi_{t}=\frac{p+q}{2},}\vspace{2mm}\\
\displaystyle{p_{t}-\lambda p_{x}=-\frac14(p^{2}+3pq),}\vspace{2mm}\\
\displaystyle{q_{t}+\lambda q_{x}=-\frac14(q^{2}+3pq),}\vspace{2mm}\\
\displaystyle{t=0:~~\phi=\ln u_{0}(x), ~p=p_{0}(x), ~q=0,}%
\end{array}
\right.
\end{equation}
where $p_{0}(x)$ is defined by (\ref{2.18}).

Noting the last equation in (\ref{2.19}), we observe that the above Cauchy
problem becomes
\begin{equation}
\label{4.6}\left\{
\begin{array}
[c]{l}%
\displaystyle{\phi_{t}=\frac12p,}\vspace{2mm}\\
\displaystyle{p_{t}-\lambda p_{x}=-\frac14p^{2},}%
\end{array}
\right.
\end{equation}
\begin{equation}
\label{4.7}t=0:~~\phi=\ln u_{0}(x),~p=p_{0}(x).
\end{equation}
On the other hand, by (\ref{2.18}) and (\ref{4.2}), we have
\begin{equation}
\label{4.8}p_{0}(x)=2\displaystyle\frac{u_{0}^{\prime}(x)}{u_{0}^{\frac{3}{2}%
}(x)}.
\end{equation}
Thus, passing through any fixed point $(0,\alpha)$ in the $x$-axis, we draw
the characteristic $\xi=\xi_{-}(t;0,\alpha)$ which is defined by (\ref{2.25}).
Along the characteristic $\xi=\xi_{-}(t;0,\alpha)$, it follows from the second
equation in (\ref{4.6}) that
\begin{equation}
\label{4.9}p\left(  t,\xi_{-}(t;0,\alpha)\right)  =\displaystyle\frac
{p_{0}(\alpha)}{1+\frac{1}{4}p_{0}(\alpha)t}.
\end{equation}
In particular,
\begin{equation}
\label{4.10}p\left(  t,\xi_{-}(t;0,x_{0})\right)  =\displaystyle\frac
{p_{0}(x_{0})}{1+\frac{1}{4}p_{0}(x_{0})t}%
\end{equation}
for $t\in\left[  0,-\displaystyle\frac{4}{p_{0}(x_{0})}\right)  $. Here we
have made use of (\ref{4.1}). By (\ref{4.1}), we find $p_{0}(x_{0})<0$, and
then $-\displaystyle\frac{4}{p_{0}(x_{0})}>0$. (\ref{4.10}) implies that
\begin{equation}
\label{4.11}p\left(  t,\xi_{-}(t;0,x_{0})\right)  \searrow-\infty
\quad\text{as}\quad t\nearrow-\displaystyle\frac{4}{p_{0}(x_{0})}.
\end{equation}
(\ref{4.11}) shows that the smooth solution of the Cauchy problem
(\ref{4.6})-(\ref{4.7}) exists only in finite time. For any $\alpha
\in\mathbb{R}$, (\ref{4.9}) always holds. Combining this and (\ref{4.11})
gives (\ref{4.3}).

We next prove (\ref{4.4}). For simplicity, we assume that
\begin{equation}
\label{4.12}p_{0}(x_{0})=\inf_{x\in\mathbb{R}}{p_{0}(x)}<0.
\end{equation}
By (\ref{4.9}), we have
\begin{equation}
\label{4.13}|p(t,x)|=\displaystyle\frac{|p_{0}(x)|_{C^{0}}}%
{1-\frac{1}{4}|p_{0}(x_{0})|t}, \quad\forall\; (t,x)\in[0,\widetilde{T}_{\max
})\times\mathbb{R}.
\end{equation}
On the other hand, it follows from the first equation in (\ref{4.6}) that
\begin{equation}
\label{4.14}\varphi(t,x)=\varphi_{0}(x)+\displaystyle\frac{1}{2}\int_{0}%
^{t}{p(\tau,x)}d\tau,
\end{equation}
and then
\begin{equation}
\label{4.15}\begin{aligned} |\varphi(t,x)|
&\leqslant|\ln{M}|+\displaystyle\frac{1}{2}
\int_0^t{\frac{|p_0(x)|_{C^0}}{1-\frac{1}{4}|p_0(x_0)|\tau}}d\tau\\
&=|\ln{M}|+
\displaystyle\frac{2|p_0(x)|_{C^0}}{|p_0(x_0)|}\ln{\left(1-\frac{1}{4}|p_0(x_0)|t\right)^{-1}},
\quad \forall\; (t,x)\in[0,\widetilde{T}_{\max})\times\mathbb{R}.
\end{aligned}
\end{equation}

We next estimate $r$ and $s$. Noting the second equation in
(\ref{2.30}) and the fact $q\equiv0$, we have
\begin{equation}
\label{4.16}s\equiv0,\quad\forall\; (t,x)\in[0,\widetilde{T}_{\max}%
)\times\mathbb{R}.
\end{equation}
Thus the Cauchy problem (\ref{2.31})-(\ref{2.32}) reduces to
\begin{equation}
\label{4.17}%
\begin{cases}
r_{t}-\lambda r_{x}=-\displaystyle\frac{3}{4}pr,\\
t=0: r=r_{0}(x)\triangleq p_{0}^{\prime}(x).
\end{cases}
\end{equation}
Thus, along the characteristic $\xi=\xi_{-}(t;0,\alpha)$, it holds that
\begin{equation}
\label{4.18}r=p_{0}^{\prime}(\alpha)\exp\left\{  -\displaystyle\frac{3}{4}%
\int_{0}^{t}{p(\tau,\xi_{-}(\tau;0,\alpha))}d\tau\right\}  .
\end{equation}
By (\ref{4.13}), we have
\begin{equation}
\label{4.19}|r(t,x)|\leq|p_{0}^{\prime}(x)|_{C^{0}} \exp\left\{
\displaystyle\frac{3|p_{0}(x)|_{C^{0}}}{|p_{0}(x_{0})|} \ln{\left(
1-\frac{1}{4}|p_{0}(x_{0})|t\right)  ^{-1}}\right\}  ,\quad\forall\;
(t,x)\in[0,\widetilde{T}_{\max})\times\mathbb{R}.
\end{equation}
In particular,
\begin{equation}
\label{4.20}r\left(  t,\xi_{-}(t;0,x_{0})\right)  =p_{0}^{\prime}(x_{0})
\exp\left\{  -\displaystyle\frac{3}{4}\int_{0}^{t}{p(\tau,\xi_{-}(\tau
;0,x_{0}))}d\tau\right\}  , \quad\forall\; t<\widetilde{T}_{\max}.
\end{equation}
Noting (\ref{4.12}), we have
\begin{equation}
\label{4.21}p_{0}^{\prime}(x_{0})=0,
\end{equation}
and then,
\begin{equation}
\label{4.22}r\left(  t,\xi_{-}(t;0,x_{0})\right)  \equiv0,\quad\forall\;
t\in[0,\widetilde{T}_{\max}).
\end{equation}

We finally estimate $R$. By (\ref{3.26}) and (\ref{4.22}), we have
\begin{equation}
\label{4.23}\begin{aligned}R\left(t,\xi_-(t;0,x_0)\right)& =\displaystyle\frac{1}{2}\left\{r\exp\{-\frac{\varphi}{2}\}+\frac{1}{4}p^2\right\}\left(t,\xi_-(t;0,x_0)\right)\\ &=\frac{1}{8}p^2\left(t,\xi_-(t;0,x_0)\right). \end{aligned}
\end{equation}
Noting (\ref{4.10}), we obtain
\begin{equation}
\label{4.24}R\left(  t,\xi_{-}(t;0,x_{0})\right)  =\displaystyle\frac{1}%
{8}\frac{p_{0}^{2}(x_{0})}{\left(  1-\frac{1}{4}|p_{0}(x_{0})|t\right)  ^{2}}
\nearrow+\infty\quad\text{as}\quad t\nearrow\widetilde{T}_{\max}.
\end{equation}

Denote
\begin{equation}
\label{4.25}x_{\ast}=\xi_{-}(\widetilde{T}_{\max};0,x_{0})\triangleq
\lim_{t\rightarrow\widetilde{T}_{\max}}\xi_{-}(t;0,x_{0}).
\end{equation}
It is easy to see that the desired (\ref{4.4}) comes from (\ref{4.24})
directly. This proves Theorem 4.1. $\qquad\quad\blacksquare$

\vskip 3mm

Theorem 4.1 is nothing but Theorem 1.3. Thus Theorem 1.3 has been
proved.

\begin{Remark}
Under the assumptions of Theorem 4.1, the scalar curvature $R$ goes to
positive infinity at algebraic rate $\displaystyle\frac{1}{(\widetilde
{T}_{max}-t)^{2}}$ as $(t,x)$ tends to the blow up point $(\widetilde{T}%
_{max},x_{\ast})$.$\qquad\quad\square$
\end{Remark}

\section{Radial solutions to hyperbolic geometric flow}

In this section, we shall investigate the radial hyperbolic
geometric flow, i.e., we shall consider the radial solution
$u=u(t,r)$ of the nonlinear wave equation (\ref{2.5}). More
precisely speaking, we consider the Cauchy problem
\begin{equation}\label{51}\left\{\begin{array}{l}u_{tt}-\triangle\ln u=0,\vspace{2mm}\\
t=0:\quad u=u_0(r),~~u_t=u_1(r),\end{array}\right.\end{equation}
where $u_0,~u_1$ are smooth functions of $r$, in which
\begin{equation}\label{52}r=\sqrt{x^2+y^2}\geqslant0.\end{equation}

In what follows, we will state and prove some formulas which are
useful in the study of radial solutions to the hyperbolic geometric
flow.

As in (\ref{2.8}), let
\begin{equation}\label{53}\psi(t,r)=\ln u(t,r).\end{equation}
Then (\ref{2.5}) becomes
\begin{equation}\label{54}\psi_{tt}-e^{-\psi}\left(\psi_{rr}+\dfrac1r\psi_r\right)+\psi_t^2=0.\end{equation}
Similar to (\ref{2.11}), we denote
\begin{equation}\label{55}v=\psi_t,\quad w=\psi_r.\end{equation}
Then, the equation (\ref{54}) can be equivalently rewritten as
\begin{equation}\label{56}\left\{\begin{array}{l}w_t-v_r=0,\vspace{2mm}\\
v_t-e^{-\psi}\left(w_r+\dfrac1r
w\right)+v^2=0.\end{array}\right.\end{equation}

Introduce
\begin{equation}\label{57}\mu=v+e^{-\frac{\psi}{2}}w,\quad \nu=v-e^{-\frac{\psi}{2}}w,\end{equation}
we have
\begin{Lemma}$\mu$ and $\nu$ satisfy the following equations
\begin{equation}\label{58}\left\{\begin{array}{l}\mu_t-\chi\mu_r=-\left(\dfrac{\mu+\nu}{2}\right)^2+\dfrac12
\left(\dfrac{\chi}{r}-\dfrac{\nu}{2}\right)(\mu-\nu),\vspace{3mm}\\
\nu_t+\chi\nu_r=-\left(\dfrac{\mu+\nu}{2}\right)^2+\dfrac12
\left(\dfrac{\chi}{r}+\dfrac{\mu}{2}\right)(\mu-\nu),\end{array}\right.\end{equation}
where
\begin{equation}\label{59}\chi=e^{-\frac{\psi}{2}}.\quad\quad\quad
\square\end{equation}\end{Lemma}

\noindent{\bf Proof.}  We calculate
\begin{equation}\label{510}\begin{aligned}\mu_t-\chi\mu_r&=\left(v+e^{-\frac{\psi}{2}}w\right)_t-
e^{-\frac{\psi}{2}}\left(v+e^{-\frac{\psi}{2}}w\right)_r\vspace{3mm}\\
&=v_t+e^{-\frac{\psi}{2}}w_t-\frac12e^{-\frac{\psi}{2}}w\psi_t-e^{-\frac{\psi}{2}}\left(v_r+
e^{-\frac{\psi}{2}}w_r-\frac12e^{-\frac{\psi}{2}}w\psi_r\right)\vspace{3mm}\\
&=v_t-\frac12e^{-\frac{\psi}{2}}w\psi_t-e^{-\psi}w_r+\frac12e^{-\psi}w\psi_r\vspace{3mm}\\
&=e^{-\psi}\left(w_r+\frac1r
w\right)-v^2-\frac12e^{-\frac{\psi}{2}}w\psi_t-e^{-\psi}w_r+\frac12e^{-\psi}w\psi_r\vspace{3mm}\\
&=-v^2+\left(\frac1r
e^{-\psi}-\frac12ve^{-\frac{\psi}{2}}+\frac12e^{-\psi}w\right)w.\end{aligned}\end{equation}
In (\ref{510}), we have made use of (\ref{56}). Noting (\ref{57}),
we obtain from (\ref{510}) that
\begin{equation}\label{511}\begin{aligned}\mu_t-\chi\mu_r&=
-\left(\frac{\mu+\nu}{2}\right)^2+\left(\frac1re^{-\frac{\psi}{2}}-\frac{\mu+\nu}{4}+\frac{\mu-\nu}{4}\right)
\frac{\mu-\nu}{2}
\vspace{2mm}\\
&=-\left(\frac{\mu+\nu}{2}\right)^2+\left(\frac{\chi}{r}-\frac{\nu}{2}\right)\frac{\mu-\nu}{2}.\end{aligned}\end{equation}

Similarly, we can prove
\begin{equation}\label{512}
\nu_t+\chi\nu_r=-\left(\frac{\mu+\nu}{2}\right)^2+\left(\frac{\chi}{r}+\frac{\mu}{2}\right)\frac{\mu-\nu}{2}\end{equation}
This proves Lemma 5.1.$\quad\quad\quad\blacksquare$

\vskip 4mm

It is easy to see that the $C^2$ solution of (\ref{2.5}) is
equivalent to the corresponding $C^1$ solution of the following
system of first order
\begin{equation}\label{513}\left\{\begin{array}{l}\psi_t=\dfrac{\mu+\nu}{2},\vspace{2mm}\\
\mu_t-\chi\mu_r=-\left(\dfrac{\mu+\nu}{2}\right)^2+\left(\dfrac{\chi}{r}-\dfrac{\nu}{2}\right)\dfrac{\mu-\nu}{2},\vspace{2mm}\\
\nu_t+\chi\nu_r=-\left(\dfrac{\mu+\nu}{2}\right)^2+\left(\dfrac{\chi}{r}+\dfrac{\mu}{2}\right)\dfrac{\mu-\nu}{2},\end{array}\right.\end{equation}
where $\chi$ is given by (5.9).

Similar to Lemma 2.2, we can prove

\begin{Lemma}It holds that
\begin{equation}\label{514}\left\{\begin{array}{l}\mu_t-(\chi\mu)_r=-\mu\nu+\dfrac{\chi(\mu-\nu)}{2r},\vspace{2mm}\\
\nu_t+(\chi\nu)_r=-\mu\nu+\dfrac{\chi(\mu-\nu)}{2r}.\end{array}\quad\quad\quad
\square\right.\end{equation}\end{Lemma}

Similar to (\ref{2.30}), let
\begin{equation}\label{515}\eta=\mu_r,\quad\gamma=\nu_r.\end{equation}
By a direct calculation, we have

\begin{Lemma}$\eta$ and $\gamma$ satisfy
\begin{equation}\label{516}\left\{\begin{array}{l}\eta_t-\chi\eta_r=-\dfrac{3\mu+\nu}{4}(\eta+\gamma)-\left(\dfrac{\mu-\nu}{4r}
+\dfrac{\chi}{r^2}\right)\dfrac{\mu-\nu}{2}+\left(\dfrac{\chi}{r}-\dfrac{\nu}{2}\right)\dfrac{\eta-\gamma}{2},
\vspace{3mm}\\
\gamma_t+\chi\gamma_r=-\dfrac{\mu+3\nu}{4}(\eta+\gamma)-\left(\dfrac{\mu-\nu}{4r}+\dfrac{\chi}{r^2}\right)\dfrac{\mu-\nu}{2}
+\left(\dfrac{\chi}{r}-\dfrac{\mu}{2}\right)\dfrac{\eta-\gamma}{2}.\end{array}\quad\square
\right.\end{equation}
\end{Lemma}

Lemmas 5.1--5.3 play an important role in the study on the global
existence and blowup phenomenon of the radial solutions of the
hyperbolic geometric flow. The key point is to establish the
estimates on the terms with the factor $\dfrac1r$ in the equations
(\ref{58}) and (\ref{516}). This kind of estimates is a difficult
and key point in the study on the theory of hyperbolic partial
differential equations, and still remains open.

\section{Some geometric properties on general Riemann surfaces}

In this section, we consider the geometric properties of solutions
of the hyperbolic geometric flow on Riemann surface of general
initial value. We first give some special explicit solutions of the
reduced equation of the hyperbolic geometric flow to highlight some
special forms of the solutions. Then we derive the reduced equation
of the hyperbolic geometric flow in more general conformal class and
give some results about certain important quantities like the total
mass or volume function which illustrate some geometric obstructions
to the long time existence of the reduced problems, as well as the
existence of periodic solutions.

Let us consider a simple case, let $\left(\mathscr{M},g\right)  $ be
a complete Riemann surface with a metric of conformal type%
\[
g=\rho\left(  x,y\right)  ^{2}\left(  dx^{2}+dy^{2}\right)  .
\]
Here we can assume that $\mathscr{M}$ is globally conformal to
$\mathbb{R}^{2}$ or $\mathbb{T}^{2}$. If we look for the solutions
of the hyperbolic geometric flow (\ref{1.1}) or (\ref{1.6}) in this
conformal class, as in Section 1, we can reduce this
system explicitly to%
\begin{equation}
\frac{\partial^{2}u}{\partial t^{2}}-\triangle\ln u=0, \label{5.2}%
\end{equation}
where $u=\rho^{2}$ and $\bigtriangleup$ is the standard Laplacian.
This can be done by a simple observation that the Ricci curvature
can be written in terms of the Gaussian curvature as $Ric\left(
g_{t}\right)  =Kg_{t}$, and the Gaussian curvature can be written as
$K=-\frac{1}{\rho^{2}}\triangle\ln\rho$.

Setting $w=\ln u$, same as (2.9) we have
\begin{equation}
\frac{\partial^{2}w}{\partial t^{2}}+\left\vert \frac{\partial
w}{\partial
t}\right\vert ^{2}-e^{-w}\triangle w=0. \label{5.3}%
\end{equation}
The equation (\ref{5.3}) is reminiscent of the traditional wave map
and semilinear hyperbolic equation in two dimension. In general, it
is not easy to handle globally. We first search for some special
solutions.

\begin{Example}
(\textbf{Solutions of separation variables}) {\em  We look for
solutions of (\ref{5.3})  with the following form
\[
w=f\left(  t\right)  +g\left(  x,y\right) .
\]
We have%
\[
\frac{d^{2}f}{dt^{2}}+\left\vert \frac{df}{dt}\right\vert ^{2}-e^{-f\left(
t\right)  }e^{-g\left(  x,y\right)  }\bigtriangleup g=0,
\]
this means%
\[
e^{f\left(  t\right)  }\left(  \frac{d^{2}f}{dt^{2}}+\left\vert \frac{df}%
{dt}\right\vert ^{2}\right)  =e^{-g\left(  x,y\right)  }\bigtriangleup g=c,
\]
where $c$ is a constant independent of $t,x,y$. So we can solve
these
two equations by%
\[
e^{f\left(  t\right)  }=\frac{c}{2}t^{2}+at+b,
\]
and%
\begin{equation}
\bigtriangleup g=ce^{g\left(  x,y\right)  },\label{5.4}%
\end{equation}
where $a,b$ are arbitrary real numbers, and $g\left(  x,y\right)  $
is a solution of (\ref{5.4}) which is the famous Liouville equation
which has arisen in complex analysis and differential geometry on
Riemann surfaces, in particular, in the problem of prescribing
curvature. We just notice that when the Riemann surface is compact,
there is no solution to (\ref{5.4}) unless
$c=0$.$\qquad\quad\square$}
\end{Example}

\begin{Example}
(\textbf{Solutions of traveling wave type}) {\em We look for
solutions of
(\ref{5.3}) of the type%
\[
w=f\left(  x-at\right)  ,
\]
where $a$ is a real number, then we have the following ordinary
differential
equation%
\[
e^{f}\left(  a^{2}f^{\prime\prime}+a^{2}f^{\prime2}\right)  -f^{\prime\prime
}=0,
\]
which can be solved explicitly in the following implicit form%
\[
a^{2}e^{f}-f=c_{1}x+c_{2}.\qquad\quad\square
\]}

\end{Example}

Let $\left(\mathscr{M},g_{0}\right)  $ be a Riemann surface with
metric $g_{0}$. We shall call metrics $g_{0}$ and $g$ point-wise
conformally equivalent if $g=ug_{0}$ for some positive function
$u\in C^{\infty}\left(\mathscr{M}\right)  $, whereas we say that
metrics $g_{0}$ and $g$ are conformally equivalent if there is a
diffeomorphism $\Phi$ of $\mathscr{M}$ and a positive function $u\in
C^{\infty}\left(\mathscr{M}\right)  $ such that $\Phi^{\ast}\left(
g\right) =ug$. Although we assume that the Riemann surface is
globally conformal to $\mathbb{R}^{2}$ in the above simple setting,
we can investigate the hyperbolic geometric flow in the same
pointwise conformal class in more general setting. Let
$\left(\mathscr{M},g_{0}\right)$ be a Riemann surface with metric
$g_{0}$, we would like to search for solutions of the hyperbolic
geometric flow in the pointwise conformal class
$\left(\mathscr{M},u\left(t, \cdot\right)  g_{0}\right)  $, where
$u\left(t, \cdot\right)  $ is a function with parameter $t$.

To make the problem more tangible, we first prove the following
lemma.

\begin{Lemma}
Let $g_{0}$ and $g$ be point-wise conformally equivalent metrics on
the Riemann surface $\mathscr{M}$ such that $g=ug_{0}$ for a
positive smooth function $u$, and $k\left(  \cdot\right)  $ is the
Gaussian curvature of $\left(\mathscr{M},g_{0}\right)$, then the
Gaussian curvature of $\left(\mathscr{M},ug_{0}\right)$ is
\begin{equation}
K\left(  x,y\right)  =\frac{1}{u}\left(  k\left(  x,y\right)  -\frac{1}%
{2}\triangle_{g_{0}}\ln u\right)  ,\label{5.5}%
\end{equation}
where $\left(  x,y\right)  \in \mathscr{M}$ is any local parameter and $\triangle_{g_{0}}%
$denotes the Laplacian with respect to metric
$g_{0}$.$\qquad\quad\square$
\end{Lemma}

\noindent\textbf{Proof.} Let $\left\{
\omega_{1}^{0},\omega_{2}^{0}\right\} $ be a local oriented
orthonormal coframe field on $\left(\mathscr{M},g_{0}\right)  $. If
we set $\omega_{i}=\sqrt{u}\omega_{i}^{0}$, then $\left\{  \omega
_{1},\omega_{2}\right\}  $ is a local oriented orthonormal coframe
field for $g$. It is worth to notice now that the Gaussian curvature
$k$ of $\left(\mathscr{M},g_{0}\right)$ is determined by the following equation%
\[
k\omega_{1}^{0}\wedge\omega_{2}^{0}=d\varphi_{12}^{0},
\]
where $\varphi_{12}^{0}$ denotes the Riemannian connection form
$g_{0}$.

Now we compute the connection $\varphi_{12}$ form of
$\left(\mathscr{M},u\left( \cdot\right)  g_{0}\right)  $. Let $du=u_{1}\omega_{1}^{0}+u_{2}\omega_{2}%
^{0}$. Then%
\begin{align*}
d\omega_{1}  &  =\left(  \frac{1}{2\sqrt{u}}du\wedge\omega_{1}^{0}-\sqrt
{u}\varphi_{12}^{0}\wedge\omega_{2}^{0}\right) \\
&  =-\left(  \frac{u_{2}}{2u}\omega_{1}^{0}-\frac{u_{1}}{2u}\omega_{2}%
^{0}+\varphi_{12}^{0}\right)  \wedge\omega_{2}.
\end{align*}
Thus
\[
\varphi_{12}=\frac{u_{2}}{2u}\omega_{1}^{0}-\frac{u_{1}}{2u}\omega_{2}%
^{0}+\varphi_{12}^{0}=\varphi_{12}^{0}-\frac{1}{2}\ast d\ln u,
\]
where $\ast$ denotes the Hodge star operator. So we conclude%
\[
K\omega_{1}\wedge\omega_{2}=d\varphi_{12}=d\varphi_{12}^{0}-\frac{1}{2}d\ast
d\ln
u=k\omega_{1}^{0}\wedge\omega_{2}^{0}-\frac{1}{2}\triangle\left(
\ln u\right)  \omega_{1}^{0}\wedge\omega_{2}^{0},
\]
and the result follows immediately. $\qquad\quad\blacksquare$

\vskip 4mm

>From (\ref{5.5}) we can reduce the hyperbolic geometric flow
(\ref{1.1}) in a
point-wise conformal class to%
\begin{equation}
\frac{\partial^{2}u}{\partial t^{2}}=\triangle_{g_{0}}\ln u-2k\left(
x,y\right)  , \label{5.6}%
\end{equation}
and moreover we can consider the Cauchy problem for (\ref{5.6})
\begin{equation}
\left\{
\begin{array}
[c]{l}
{\displaystyle \frac{\partial^{2}u}{\partial
t^{2}}=\triangle_{g_{0}}\ln u-2k\left(
x,y\right)} ,\\
t=0:\;\;u=u_0(x,y),\quad u_t=u_1(x,y),
\end{array}
\right.  \label{5.7}%
\end{equation}
where $u_{0}$ is a positive function and $u_{1}$ is an arbitrary
function on $\mathscr{M}$.

Let us now consider (\ref{5.6}) or (\ref{5.7}). As investigated by
Kong and Liu \cite{kl}, we know that the hyperbolic geometric flow
on three dimension has a deep relation with the Einstein equations.
It is a challenging problem to find periodic solutions and blowing
up phenomena of the hyperbolic geometric flow as in the Einstein
equation case. It is interesting that in two dimensional case we can
prove the following theorem which gives some obstructions to the
existence of smooth long time solutions and periodic solution of
(\ref{5.7}).

In order to state our result, we introduce
\begin{Definition}
A solution of (\ref{5.7}) is regular, if it is positive and
smooth.$\qquad\quad\square$
\end{Definition}

\begin{Theorem}
Let $\left(\mathscr{M},g_{0}\right)  $ be a compact Riemannian
surface with metric $g_{0}$, $\chi\left(\mathscr{M}\right)  $
denotes its Euler characteristic number. If $u\left(t,x,y\right) $
is a regular solution of (\ref{5.7}). We have

(a) If $\chi\left(\mathscr{M}\right)  >0$, then any solution of the
reduced problem (\ref{5.7}) must blow up in finite time for any
initial value $(u_{0},u_{1})$;

(b) If $\chi\left(\mathscr{M}\right)  \neq0$, there is no periodic
regular solution of equation (\ref{5.6}). Moreover in the case
$\chi\left(\mathscr{M}\right)=0$, if there is a periodic solution
$u\left(t,x,y\right)$ of (\ref{5.6}), then we must have
$\int\limits_{\mathscr{M}}u\left(t,x,y\right)dV_{g_{0}}=c$ for some
positive constant $c$;

(c) If $\chi\left(\mathscr{M}\right)=0$, and $u\left(t,x,y\right) $
is a regular solution of (\ref{5.7}). Then if
$\int\limits_{\mathscr{M}}u_{1}\left( x,y\right) dV_{g_{0}}<0$,
$u\left(t,x,y\right)  $ must blow up in finite time.
$\qquad\quad\square$
\end{Theorem}

\noindent\textbf{Proof.} Let $u\left(t,x,y\right)  $ be a regular
solution of (\ref{5.6}) or (\ref{5.7}). We denote the volume of
$\mathscr{M}$ with respect to the metric $u\left(t,x,y\right)
g_{0}$ by
\[
V\left(  t\right)  =%
{\displaystyle\int\limits_{\mathscr{M}}}
u\left(t,x,y\right)  dV_{g_{0}}.
\]
Then taking integration on both sides of (\ref{5.6}) and using Gauss-Bonnet
formula, we have%
\[
\frac{d^{2}V\left(  t\right)
}{dt^{2}}=-4\pi\chi\left(\mathscr{M}\right)  ,
\]
this means%
\begin{equation}
V\left(t\right)  =-2\pi\chi\left(\mathscr{M}\right)  t^{2}+c_{1}t+c_{2} \label{5.8}%
\end{equation}
for some constants $c_{1}$ and $c_{2}$. If we consider the Cauchy
problem
(\ref{5.7}), we can easily get the numbers $c_{1}$ and $c_{2}$, $c_{2}=%
{\displaystyle\int\limits_{\mathscr{M}}}
u_{0}\left(  x,y\right)  dV_{g_{0}}$ and $c_{1}=%
{\displaystyle\int\limits_{\mathscr{M}}}
u_{1}\left(  x,y\right)  dV_{g_{0}}$. From (\ref{5.8}) we conclude
(a) by the assumption $\chi\left(\mathscr{M}\right)  >0$. Otherwise,
for some sufficiently large number $T$, we have $V\left(  T\right)
\leq0$, which is impossible.

To prove (b), we just notice that if $u\left(t,x,y\right)  $ is a
periodic solution of (\ref{5.6}), then $V\left( t\right)  $ must be
a periodic function of parameter $t$, which is impossible unless
$\chi\left(\mathscr{M}\right) =0$ and $c_{1}=0$ in (\ref{5.8}).

The conclusion (c) can also be concluded from the expression
(\ref{5.8}). In this case we have $$V\left(  t\right)
=c_{1}t+c_{2},$$ but
\[
c_{1}=%
{\displaystyle\int\limits_{\mathscr{M}}}
u_{1}\left(  x,y\right)  dV_{g_{0}}<0,
\]
so $u\left(t,x,y\right)$ must blow up in finite time. This proves
Theorem 6.1. $\qquad\quad\blacksquare$

\vskip 4mm

>From the above theorem, we see that problem (\ref{5.7}) or equation
(\ref{5.6}) have essential geometric obstruction to the existence of
long time solution or periodic solution. However, as in the study of
Ricci flow on Riemann surfaces, we consider the following normalized
equation for the equation
(\ref{5.6})%
\begin{equation}
\left\{
\begin{array}
[c]{l}%
{\displaystyle\frac{\partial^{2}u}{\partial
t^{2}}=\triangle_{g_{0}}\ln u-2k\left(
x,y\right)  +\frac{4\pi\chi\left(\mathscr{M}\right)  }{V\left(  g_{0}\right)  },}\\
t=0:\;\;u=u_0(x,y),\quad u_t=u_1(x,y),
\end{array}
\right.  \label{5.9}%
\end{equation}
where $V\left(  g_{0}\right)  $ denotes the volume of $\mathscr{M}$
with respect to metric $g_{0}$. We have

\begin{Theorem}
Let $\left(\mathscr{M},g_{0}\right)  $ be a compact Riemann surface
with metric $g_{0}$, $\chi\left(\mathscr{M}\right)  $ denotes its
Euler characteristic number. If $u\left(t,x,y\right) $ is a regular
solution of (\ref{5.9}), then we have

(i) If $%
{\displaystyle\int\limits_{\mathscr{M}}}
u_{1}\left(  x,y\right)  dV_{g_{0}}<0$, then the volume function
$V\left( t\right)  $ is linearly contracting and
$u\left(t,x,y\right) $ must blow up in a finite time;

(ii) If $%
{\displaystyle\int\limits_{\mathscr{M}}}
u_{1}\left(  x,y\right)  dV_{g_{0}}>0$, then the volume function $V\left(
t\right)  $ is linearly expanding.$\qquad\quad\square$
\end{Theorem}

This theorem can be proved in a manner similar to the proof of
(\ref{5.8}), here we omit the details.

Theorem 6.2 reveals that the hyperbolic geometric flow  has very
different features from the traditional Ricci flow, because even for
the normalized equation (\ref{5.9}) we can not get the long-time
existence of solutions which depends on the velocity $u_{1}$. This
also gives an evidence that, by choosing the initial velocity of the
hyperbolic geometric flow, we may obtain the long time existence of
the corresponding solution. We believe that the above results may
reveal some interesting features of the Einstein equations about the
expansion of the Universe, if one could connect the foliation of
solution of the Einstein equations in three dimension with the
solution of (\ref{5.7}) or (\ref{5.9}). We will pursue this problem
in a forthcoming paper.

\vskip 5mm

\noindent{\Large \textbf{Acknowledgements.}} The work of Kong was supported in
part by the NNSF of China (Grant No. 10671124) and the NCET of China (Grant
No. NCET-05-0390); the work of Liu was supported in part by the NSF and NSF of China.


\begin{thebibliography}{9}                                                                                                %

\bibitem {dkl}Wen-Rong Dai, De-Xing Kong and Kefeng Liu, \textit{Hyperbolic
geometric flow (I): short-time existence and nonlinear stability},
to appear in Pure and Applied Mathematics Quarterly: special issue
dedicated to M. Atiyah's 80th birthday.

\bibitem {dkl2}Wen-Rong Dai, De-Xing Kong and Kefeng Liu, \textit{Dissipative hyperbolic
geometric flow}, arXiv:0709.2542.

\bibitem{h} Jia-Xing Hong, \textit{The global smooth solutions of Cauchy
problems for hyperbolic equation of Monge-Ampere type}, Nonlinear
Analysis, Theory, Method \& Applications \textbf{24} (1995),
1649-1663.

\bibitem {k}De-Xing Kong, \textit{Maximum principle in nonlinear hyperbolic
systems and its applications}, Nonlinear Analysis, Theory, Method \&
Applications \textbf{32} (1998), 871-880.

\bibitem {k2} De-Xing Kong, \textit{Hyperbolic geometric flow}, the
Proceedings of ICCM 2007, Vol. I\!I, Higher Educationial Press,
Beijing, 2007, 95-110.

\bibitem {kl}De-Xing Kong and Kefeng Liu, \textit{Wave character of metrics
and hyperbolic geometric flow}, J. Math. Phys. 48 (2007),
103508-1-103508-14.

\bibitem {p} R. Penrose, \textit{Gravitational collapse and space-time
singularities}, Phys. Rev. Lett. \textbf{14} (1965), 57-59.

\bibitem {y}R. Schoen and S.-T. Yau, Lectures on Differential Geometry,
International Press, Cambridge, MA, 1994.
\end{thebibliography}
\end{document}